\documentclass[a4paper,12pt,reqno]{amsart}
\usepackage[utf8]{inputenc}
\usepackage{amssymb,amsthm}
\usepackage{enumerate}
\usepackage{hyperref}
\usepackage[margin=3cm]{geometry}
\usepackage{graphicx}
\usepackage{color}

\renewcommand{\geq}{\geqslant}
\renewcommand{\epsilon}{\varepsilon}
\renewcommand{\leq}{\leqslant}

\newtheorem{theorem}{Theorem}[section]

\newcommand{\R}{\mathbb{R}}
\newcommand{\C}{\mathbb{C}}

\DeclareMathOperator{\divv}{div}

\title[A multidimensional fundamental theorem of calculus]{An elementary multidimensional fundamental theorem of calculus}

\author[Joaquim Bruna]{Joaquim Bruna}
\address{Departament de Matem\`{a}tiques, Edifici Cc, Universitat Aut\`{o}noma de Barcelona, 08193 Cerdanyola del Vall\`es, Barcelona, Spain.}
 \email{joaquim.bruna@uab.cat}

\date{}

\begin{document}

\begin{abstract} We discuss a  version of the fundamental theorem of calculus in several variables and some applications, of potential interest as a teaching material in  undergraduate courses.
\end{abstract}

\maketitle

\section{Introduction and main result}

In standard undergraduate courses on one variable calculus, differentiation and integration are presented as inverse processes, as stated in the \emph{fundamental theorem of calculus}: if  $F$ is differentiable at every point $x\in [a,b]$ and $F'=f$ is integrable in $[a,b]$, then 
$$
F(x)-F(a)=\int_a^x f(t)\, dt,\quad  a\leq x\leq b.
$$
This holds both in the context of Riemann and Lebesgue integration (see \cite{rud} for a proof in the Lebesgue integration context).

\vskip 0.3 cm
In this note we provide a multidimensional version of this statement. The proof is straightforward and can  be included in an undergraduate course of multidimensional  calculus.   

\vskip 0.3 cm
The first basic concept is that of \emph{interval or cube function} $\Phi$ on a domain  $U\subset \R^n$. By an interval $Q$ in $U$ we understand a set $Q=\prod_{j=1}^n I_j$ with the $I_j=I$, $I$ a  one-dimensional closed interval, that is the faces of $Q$ are parallel to the coordinate axis.  An interval function  is 
a map defined on  all intervals $Q\subset \R^n$ assigning to each $Q$ a real or complex number $\Phi(Q)$ with the property that
$$\Phi (Q)=\sum_i \Phi(Q_i),$$
whenever $(Q_i)$ is a finite partition of $Q$, that is, $Q=\cup Q_i$ and the $Q_i$ have disjoint interiors.  The easiest example is
$$\Phi_f(Q)=\int_Q f\, dx,$$
where $f$ is locally integrable.  

We have in mind two other examples. For the first one, assume that $T:U\rightarrow V$ is a measurable homeomorphism between two domains such that $m(T(A))=0$ if $m(A)=0$. Here and in the following $m(A)$ denotes the Lebesgue measure of $A$. Then $\Phi(Q)=m(T(Q))$ is an interval function, because images of the faces have zero measure. For the second one assume that $F$ is a continuous vector field in the plane or in space and 
set
$$\Phi(Q)=\int_{\partial Q}\langle F,N\rangle \, dm_{n-1},$$
the flow of $F$ through the boundary $\partial Q$ oriented with the outward normal $N$. Then $\Phi(Q)$ is an interval function. This is because  if $Q_i,Q_j$ are two intervals with a  face $S$ in common, the outward normals are opposite one each other. 

Notice that a Dirac delta at a point $a\in U$, that is, $\Phi(Q)=1$ if $a\in Q$ and zero otherwise, is not an interval function according to our definition, because if $a$ is a boundary point of both $Q_1,Q_2$ then $\Phi(Q_1)=\Phi(Q_2)=\Phi(Q)=1$.

In dimension $n=1$, with $U=(a,b)$, if $g$ is defined on $(a,b)$, it is immediately seen that 
\begin{equation}\label{dim1}
\Phi([c,d])=g(d)-g(c),
\end{equation}
defines an interval function on $(a,b)$. Indeed, a decomposition of $[c,d]$ into pieces $Q_i$ amounts to a selection of intermediate points (the end-points of the $Q_i$) $c=t_0<t_1<\cdots<t_m=d$, and then
 $$\Phi(I)=g(d)-g(c)=\sum_i g(t_{i+1})-g(t_i)=\sum_i \Phi(Q_i).$$
 Conversely, given an interval function defined on $(a,b)$ and $p\in (a,b)$, the function
 $$g(x)=\Phi([p,x]), p\leq x; g(x)=-\Phi([x,p]), x\leq p,$$
 satisfies \eqref{dim1}. Thus there is a one-to-one correspondence between interval functions and classical functions.

The second basic concept is that of density. For an interval function $\Phi$ we define its \emph{upper density}
$$
\overline{D}_{\Phi}(x)=\limsup_{x\in Q, \delta(Q)\to 0} \frac{\Phi(Q)}{m(Q)}=
\inf_{\varepsilon}\sup_{\delta(Q)\leq \varepsilon, x\in Q}\frac{\Phi(Q)}{m(Q)}
,$$
where $m(Q)$ denotes the measure of $Q$ and $\delta(Q)$ its diameter.
 Analogously the lower density is defined
$$\underline{D}_{\Phi}(x)=\liminf_{x\in Q, \delta(Q)\to 0} \frac{\Phi(Q)}{m(Q)}=
\sup_{\varepsilon}\inf_{\delta(Q)\leq \varepsilon, x\in Q}\frac{\Phi(Q)}{m(Q)}.$$
In case both are finite and equal we say that $\Phi$ has a \emph{finite density} $D_{\Phi}(x)$ at $x$. 

For instance, if $f$ is continuous in $U$, the density of $\Phi_f$ is $f$ at all points. Indeed, given $\varepsilon>0$ there is $\tau$ such that $|f(y)-f(x)|\leq\varepsilon$ if $|x-y|\leq\tau$. Then, if $\delta(Q)<\tau, x\in Q$ one has $|f(y)-f(x)|\leq\varepsilon$ for all $y\in Q$ so
$$|\frac{\Phi_f(Q)}{m(Q)}-f(x)|=|\frac{1}{m(Q)}\int_Q (f(y)-f(x))\,dy|\leq \frac{1}{m(Q)}\int_Q |f(y)-f(x)|\,dy\leq\varepsilon.$$
Thus
$$\lim_{\delta(Q)\to 0} \frac{\Phi_f(Q)}{m(Q)}=f(x).$$
A deeper result is Lebesgue's differentiation theorem (see \cite{loja}) stating that $\Phi_f$ has density $f(x)$ at almost all points $x\in U$ under the sole assumption that $f$ is locally integrable.

For a better understanding of the density consider the following example. Assume $U=(0,1)\times (0,1)$ and let $L=\{(x,x), 0<x<1\}$ be the diagonal. Define $\Phi(Q)$ as the length of $L\cap Q$, clearly an interval function. Then $D_{\Phi}(x)=0$ for $x\notin L$ while $\underline{D}_{\Phi}(x)=0, \overline{D}_{\Phi}(x)=+\infty$ for $x\in L$.

In dimension one, if $\Phi$ is given by $\eqref{dim1}$, $\Phi$ has a finite density at $x$ if and only $g$ is differentiable at $x$, because if $x\in [c,d]$
$$\frac{g(d)-g(c)}{d-c}=\frac{g(d)-g(x)}{d-x}\frac{d-x}{d-c}+\frac{g(x)-g(c)}{x-c}\frac{x-c}{d-c}.$$

Next elementary theorem seems to be unnoticed, to the best of author's knowledge. It holds both in the context of Riemann and Lebesgue's integration.
\vskip 0.3 cm

\begin{theorem}\label{ftoc}  If an interval function $\Phi$ has a finite upper density $\overline{D}_{\Phi}$  {\bf at every point} and $\overline{D}_{\Phi}$  is locally integrable, then for every cube $Q\subset U$
$$\Phi(Q)\leq  \int_Q \overline{D}_{\Phi}(x) \, dx.$$
Analogously, if $\Phi$ has a finite lower density $\underline{D}_{\Phi}$  {\bf at every point} and $\underline{D}_{\Phi}$  is locally integrable, then
$$\Phi(Q)\geq  \int_Q \underline{D}_{\Phi}(x) \, dx.$$
Thus,
$$\Phi(Q)= \int_Q D_{\Phi}(x) \, dx,$$
whenever $\Phi$ has a finite integrable density at every point.
\end{theorem}

In an informal way, if $\Phi(Q)$ is of the order of $f(x) m(Q)$ for infinitesimal cubes $x\in Q$, then $\Phi(Q)= \int_Q f(x)\, dx$ for big cubes.

In dimension one, in view of the remark before the theorem, this is the fundamental theorem of calculus stated in the beginning.

As a corollary we may state:

{\bf Corollary.} {\it For an interval function $\Phi$ and a continuous function $f$ on $U$ the following two statements are equivalent}:

$$\lim_{x\in Q,\delta(Q)\to 0}\frac{\Phi(Q)}{m(Q)}=f(x), \quad \Phi(Q)=\int_Q f(x)\, dx.$$
\vspace{0,5cm}

We point out some remarks. First, it is essential, as in one variable, that the density is assumed to exist at {\bf every} point. If it exists just a.e. then the theorem does not hold. Secondly, in other type of results the a.e. existence of the density is actually {\bf proved} like in  Lebesgue's differentiation theorem quoted before. In fact, the interval functions $\Phi_f$ are characterized as those being \emph{absolutely continuous}, meaning that for every $\varepsilon>0$ there exists $\delta>0$ such that $\sum_i |\Phi(Q_i)|<\varepsilon$ whenever $Q_i$ are non-overlapping cubes and $\sum_i m(Q_i)<\delta$. So the result can be rephrased by saying that interval functions having finite integrable density at all points are automatically absolutely continuous. A reference for all these results is \cite{loja}.

\section{Proof}

Let $Q\subset U$ be a cube and let us break it into $2^n$ cubes $S_i$ of equal measure. Since $\Phi(Q)=\sum \Phi(S_i)$, one has
$$\Phi(S_i)\geq \frac{\Phi(Q)}{2^n},$$
whence 
$$\frac{\Phi(S_i)}{m(S_i)}\geq \frac{\Phi(Q)}{m(Q)},$$
for at least one $i$. Repeating the argument  we find  a sequence $Q_k$ of cubes, $Q_k\subset Q$, shrinking to some point $p\in Q$ such that
$$\frac{\Phi(Q_k)}{m(Q_k)}\geq \frac{\Phi(Q)}{m(Q)}.$$
Therefore $\overline{D}_{\Phi}(p)\geq \frac{\Phi(Q)}{m(Q)}$.
So $\Phi(Q)\leq  \overline{D}_{\Phi}(p)\, m(Q)$ for some point $p\in Q$. Similarly, $\Phi(Q)\geq  \underline{D}_{\Phi}(p)\, m(Q)$  for some point.
This holds for all cubes. So $\Phi\leq 0$ if $\overline{D}_{\Phi}\leq 0$. Similarly, $\Phi\geq 0$ if $\underline{D}_{\Phi}\geq 0$.

We complete now the proof in the Riemann integration context, where just the definition of the Riemann integral is used.

Let $(Q_i)$ be a  partition of $Q$; then
$$\Phi(Q)=\sum_i \Phi(Q_i)\leq \sum_i  \overline{D}_{\Phi}(p_i) m(Q_i).$$
Therefore, if $\overline{D}_{\Phi}$ is Riemann integrable, it follows that

$$\Phi(Q)\leq \int_Q \overline{D}_{\Phi}(x) \, dx.$$
In a similar way we see that
$$\Phi(Q)\geq \int_Q \underline{D}_{\Phi}(x)\, dx,$$
and so the theorem is proved when the density is Riemann integrable.

Assume now  that $\overline{D}_{\Phi}$  is Lebesgue integrable on $U$.  
We may assume $\Phi$ real-valued and use semi-continuous functions as in \cite{rud}. Recall that a function $g$ is called lower semi-continuous at a point $p$ if $\liminf_{x\to p}g(x)\geq g(p)$ and upper semi-continuous if  $\liminf_{x\to p}g(x)\leq g(p)$.

 Given $\varepsilon>0$, by the Vitali-Carathedory theorem, there is a lower semi-continuous function $v$ such that $\overline{D}_{\Phi}\leq v$ and $\int_U (v-\overline{D}_{\Phi})\, dx<\varepsilon$. Define
$$\Psi(Q)=\int_Q v\, dx-\Phi(Q).$$
Then, $v$ being lower semi-continuous,
$$\underline{D}_{\Psi}(x)=\liminf_{x\in Q}\frac{\Psi(Q)}{m(Q)}\geq\liminf_{x\in Q}\frac{1}{m(Q)}\int_Q v\, dx-\limsup_{x\in Q}\frac{\Phi(Q)}{m(Q)})\geq v(x)-\overline{D}_{\Phi}(x)\geq 0.$$
Therefore   $\Psi(Q)\geq 0$, whence
$$\Phi(Q)\leq \int_Q v\, dx=\int_Q \overline{D}_{\Phi}\, dx+\int_Q (v-\overline{D}_{\Phi})\, dx<\int_Q \overline{D}_{\Phi}\, dx+\varepsilon.$$
Since $\varepsilon$ is arbitrary, this shows that $\Phi(Q)\leq \int_Q \overline{D}_{\Phi}\, dx$ and applying the same argument to $-\Phi$ we are done.

\section{Applications}

1. As a first application we indicate a  simplified proof of a version of the change of variables formula, with minimal assumptions and not relying in the one-dimensional version and Fubini's theorem, the one stated in Theorem 7.26 in \cite{rud}:

\begin{theorem} Let $T:U\rightarrow V$ be an homeomorphism between two domains in $\R^n$, differentiable at every point $x\in U$. Assume that  
$|\det dT(x)|$ is  integrable on $U$. Then for a positive measurable function $f$ in $V$ one has
\begin{equation}\label{cdv}
\int_V f(y)\, dy=\int_U f(T(x)) |\det dT(x)|\, dx.
\end{equation}
\end{theorem}
Note that the assumption $|\det dT(x)|\neq 0$ is not made, so this version includes Sard's theorem.

We modify the proof in \cite{rud} replacing  the more advanced Radon-Nikodym differentiation theorem for absolutely continuous measures by theorem \ref{ftoc}.

First, Lemma 7.25 in \cite{rud} proves that $T$ maps sets of measure zero to sets of measure zero. As a consequence, 
$$\Phi(Q)=m(T(Q))$$
is an interval function. 

Secondly, theorem 7.24 in \cite{rud} proves that $\Phi$ has density $|\det dT(p)|$ at every point $p$. In fact, the proof in \cite{rud} uses balls, but it is easily checked that it holds for cubes too. 
We explain the basic idea for completeness. By hypothesis, we can approximate $T(p+h)$ near $p$ by $L=L(p+h)=T(p)+dT(p)(h)$, and use that $m(L(Q))=|\det L|\, m(Q)$ for all affine maps. 

One has

$$
T(p+h)=L+ E, |E|\leq \tau(|h|)|h|,\tau(h)\to 0.
$$
 with decreasing $\tau(t)$ as $t\to 0$.

\begin{figure} \begin{center}
  \includegraphics[width=.7\textwidth]{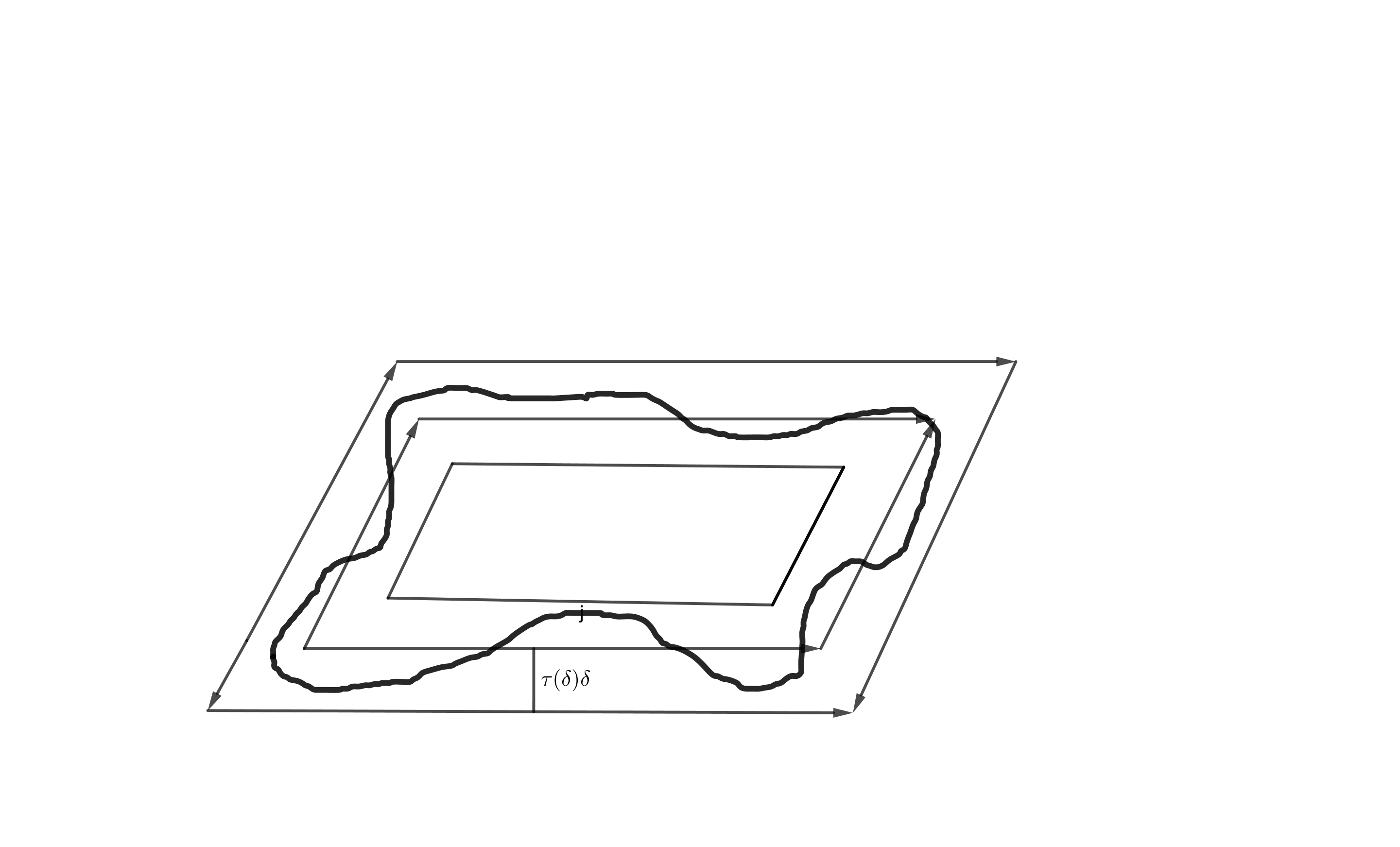}
  \caption{}\label{fig:canvivariable}
   \end{center}\end{figure}

Let $v_j=\frac{\partial T}{\partial x_j}(p), j=1,\cdots, n$ be the columns of $dT(p)$. If $Q$ has side $\delta$, $L$ maps $Q$ onto a parallelepiped $P$ with spanning vectors $\delta v_j$, whose measure is

 $$m(P)=|\det dT(p)| \delta^n= |\det dT(p)| m(Q),$$
so let us compare $T(Q)$ with $P=L(Q)$. Since $|T-L|=|E|\leq \tau(|h|)|h|$,
it is clear that $T(Q)$ is included  in a parallelepiped $P_1$ concentric with $P$ with spanning vectors $(\delta+o(\delta)) v_j$, whose measure is 
$$\delta^n |\det dT(p)|+o(\delta^n).$$
Again by $|T-L|=|E|\leq \tau(|h|)|h|$, the boundary $b(T(Q))=T(bQ)$ is at distance less than $\tau(\delta)\delta$ from $bP$. Since $T$ is an homeomorphism, this implies that  $T(Q)$ contains a parallelepiped concentric with $P$ with spanning vectors $(\delta-o(\delta)) v_j$ (see figure \ref{fig:canvivariable}  for $n=2$; a rigorous proof of this fact relies on Brouwer's fixed point theorem and can be found in lemma 7.23 of \cite{rud}), whose measure is
$$\delta^n |\det dT(p)|-o(\delta^n).$$
Altogether, since $m(Q)=\delta^n$,
$$ m(Q) |\det dT(p)|-o(m(Q))\leq m(T(Q))\leq m(Q) |\det dT(p)|+o(m(Q)),$$
proving that $\Phi$ has density $|\det dT(p)|$ at $p$.

 By theorem \ref{ftoc}, one has
 $$m(T(A))=\int_A |\det dT(x)|\, dx,$$
 when $A$ is a cube, whence when $A$ is a finite union of cubes too. Since every open set is a countable union of cubes, by the monotone convergence theorem this holds when $A$ is an open set and in turn when $A$ is a countable intersection of open sets, a $G_{\delta}$ set. Since every measurable set differs from s $G_{\delta}$ set in a set of zero measure and $T$ preserves those, we conclude that this holds for all measurable sets $A\subset U$, that is, \eqref{cdv}holds for the characteristic function of a measurable set. By linearity it then holds for simple functions, and by the monotone convergence theorem again, for a general measurable function.

 {\bf Remark.} As a first remark for the instructor, in case $|\det dT(x)|\neq 0$ for all $x\in U$, the use of Brouwer's fixed point theorem can be avoided as follows: 
 
  By the inclusion $T(Q)\subset P_1$, 
 $$m(T(Q))\leq m(Q) |\det dT(p)|+o(m(Q)),$$
 implying $\overline{D}_{\Phi}(p)\leq |\det dT(p)|$. Then theorem \ref{ftoc} implies
 $$m(T(A))\leq \int_A |\det dT(x)|\, dx,$$
 for all cubes, and as before this leads to 
 $$\int_V f(y)\, dy\leq \int_U f(T(x)) |\det dT(x)|\, dx.$$
 But since the same inequality applies to the inverse $f^{-1}$, the result follows.

{\bf Remark.} As a second remark, to be eventually combined with the previous one, a  proof in the context of Riemann integration can be further simplified as follows.  To prove \eqref{cdv} say for a continuous function $f$ with compact support, introduce 
$$\Psi(Q)=\int_{T(Q)} f\, dy.$$
The continuity of $f$ implies
$$D_{\Psi}(p)=f(T(p))D_{\Phi}(p),$$
so $D_{\Psi}(p)= f(T(p)) |\det dT(p)|$. This leads using theorem \ref{ftoc} to
$$\int_{T(Q)} f(y) \,dy= \int_Q f(T(x)) |\det dT(x)|\, dx,$$
for all cubes. If $K$ is the support of $f$, the compact $T^{-1}(K)$ can be covered by a finite number of cubes $Q$, so \eqref{cdv} follows.

2. As a second application we analyze the divergence theorem. Assume that $F$ is a continuous vector field in space and 
set
$$\Phi(Q)=\int_{\partial Q}\langle F,N\rangle \, dm_{n-1},$$
the flow of $F$ through the boundary $\partial Q$ oriented with the outward normal $N$. We mentioned before that $\Phi$ is indeed an interval function. If its density exists, we call it \emph{the divergence} $\divv F$ of $F$. If integrable, the theorem implies
$$\int_{bQ}\langle F,N\rangle \, dm_{n-1}=\int_Q \divv F\, dm_n,$$
and the same holds with $Q$ replaced by a finite union of cubes. From this it follows by approximations that the same holds with $Q$ replaced by a domain with piece-wise regular boundary (details can be found in \cite{bruna}).

 If $F$ is differentiable with components $F_i$, let us check that the density  $\divv F$ exists at every point and equals $\langle \nabla,F\rangle=\sum_i D_iF_i$. 

First  consider an affine  field $F(x)=M(X-P)$, where $M=(m_{ij})$ is  a constant matrix and $X,P$ are the column vectors $x^t, p^t$, and let us compute  the flux across the boundary $\partial Q$ of a parallelepiped in space spanned by $3$ vectors $v_1,v_2,v_3$
$$Q=\{p'+t_1 v_1+ t_2 v_2+ t_3 v_3, 0\leq t_i\leq 1 \},$$ containing $p$, oriented by the outward normal. $F$ differs from $M(X-P')$ by a constant field, which obviously has zero flux, so we can replace $p$ by $p'$ and assume $p'=0$.   On the face $t_3=1$, the basis $v_1,v_2$ is positively oriented and  the flux is
$$\int_0^1 \int_0^1 \det (M(t_1v_1+t_2v_2+v_3),v_1,v_2)\, dt_1\, dt_2,$$
while on the face $t_3=0$ it is
$$- \int_0^1 \int_0^1 \det (M(t_1v_1+t_2v_2),v_1,v_2)\, dt_1\, dt_2.$$
Therefore they add up to
$$\det (M(v_3),v_1,v_2)$$
If $M(v_3)=\sum_i \lambda_i v_i$, this equals $\lambda_3 \det (v_3,v_1,v_2)$. The same applies to the other two couples of opposite sides, whence the flux is exactly
$$\operatorname{trace}(M) \det (v_1,v_2,v_3),$$
the trace of $M$ times the volume of $Q$. 

Now let $F$ be a differentiable field at $p$, $Q$  a cube  of size $\delta$ containing $p$. As before we expand $F$ around $p$

$$F(x)=F(p)+dF(p)(X-P)+E, E=o(|x-p|).$$
The contribution to the flux of $F$ across $\partial Q$ of the constant field $F(p)$ is zero, that of the linear field $dF(p)(X-P)$ is the trace of $dF(p)$ times $m(Q)$ while that of $E$ is $o(\delta^n)$, whence the flux equals
$$(D_1F_1+\cdots+D_nF_n)(p) m(Q)+ o(m(Q)),$$
thus proving that the density is $\langle \nabla,F\rangle$.

Upon replacement of the field $F=(A,B)$ by $JF=(-B.A)$ the divergence theorem in the plane amounts to Green's formula. Using the language of line integrals, if $P\,dx+Q\,dy$ is a differentiable $1$-form  and $Q_x-P_y$ is integrable one has
$$\int_{bU} P\,dx+Q\,dy=\int_U(Q_x-P_y)\, dA,$$
with no assumption needed separately for $Q_x,P_y$. A particular case are complex line integrals 
$$\Phi(Q)=\int_{bQ} f(z)\,dz,$$
for $f$ continuous in the complex plane $\C$. For differentiable $f$ the density is $\overline{\partial}f$, and so if integrable one has
$$\int_{bU} f(z)\,dz=\int_U  \overline{\partial}f(z)\, dA(z).$$

Other general versions of Green's theorem with minimal assumptions are known, but the proofs are far from elementary (see \cite{fesq} and references herein).

3. In a surface $S$ in $\R^3$ oriented by a unit normal field $N$ one can define cubes as those which are so in a local chart. If $F$ is a continuous field, the circulation
$$\Phi(Q)=\int_{bQ}\langle F,T\rangle \, ds,$$
defines an interval function. If $F$ is differentiable, one can show along the same lines that the density is $\langle\nabla\times F,N\rangle$ and one gets Stoke's theorem with minimal assumptions (see the details in the book \cite{bruna}).

\subsection*{Acknowledgements}
The author is partially supported by the  Ministry of Science
and Innovation--Research Agency of the Spanish Government through
grant PID2021-123405NB-I00 and by AGAUR, Generalitat de Catalunya,
through  grant 2021-SGR-00087. The author thanks his colleague Juli\`{a} Cuf\'{i} for valuable comments.

\end{document}